\documentclass{article}
\begin{document}
\textbf{\Large \centerline{Generalizations of Ramanujans Continued fractions}} \vskip .10in
\[
\]
\centerline{\bf Nikos Bagis}
\centerline{bagkis@hotmail.com}
\[
\]
\begin{quote}

\centerline{\bf abstract\rm}
In this article we continue a previous work in which we have generalized the Rogers Ramanujan continued fraction (RR) introducing what we call, the Ramanujan-Quantities (RQ). We use the Mathematica package to give several modular equations for certain cases of Ramanujan Quantities-(RQ). We also give the modular equations of degree 2 and 3 for the evaluation of the first derivative of Rogers-Ramanujan continued fraction. More precicely for certain classes of (RQ)'s we show how we can find the corresponding continued fraction expansions-S, in which we are able to evaluate with numerical methods some lower degree modular equations of this fraction and its derivatives. A systematicaly method for evaluating theoriticaly certain (RQ)'s (not for all) and their derivatives, with functions used by Ramanujan himself, is  presented. We give applications and several results.          
\[
\]     
\bf keywords: \rm{Ramanujan; Continued Fractions; Quantities; Modular Equations; Derivatives; Evaluations}

\end{quote}

\section{Definitions and Introductory Results}

\label{intro}
In this article we will define and study expressions that rise from continued fractions, analogous to that of Rogers-Ramanujan (RR), Ramanujan's Cubic (RC), Ramanujan-Gollnitz-Gordon (RGG). The results are new since no work have been done in this area and most of them are experimental observations.\\ 
The focused quantities are 
\begin{equation}
R(a,b,p;q)=q^{-(a-b)/2+(a^2-b^2)/(2p)}\frac{\prod^{\infty}_{n=0}(1-q^aq^{np})(1-q^{p-a}q^{np})}{\prod^{\infty}_{n=0}(1-q^bq^{np})(1-q^{p-b}q^{np})} , 
\end{equation}
where $a$, $b$, $p$ are positive rationals such that $a+b<p$.\\
General Theorems such\\
$$
\frac{q^{B-A}}{(1-a_1b_1)+}\frac{(a_1-b_1q_1)(b_1-a_1q_1)}{(1-a_1b_1)(q_1^2+1)+}\frac{(a_1-b_1q_1^3)(b_1-a_1q_1^3)}{(1-a_1b_1)(q_1^4+1)+}\frac{(a_1-b_1q_1^5)(b_1-a_1q_1^5)}{(1-a_1b_1)(q_1^6+1)+\ldots}=
$$
\normalsize
$$
=\frac{\prod^{\infty}_{n=0}(1-q^aq^{np})(1-q^{p-a}q^{np})}{\prod^{\infty}_{n=0}(1-q^bq^{np})(1-q^{p-b}q^{np})} 
$$
where $a_1=q^A$, $b_1=q^B$, $q_1=q^{A+B}$, $a=2A+3p/4$, $b=2B+p/4$ and $p=4(A+B)$, $|q|<1$, are proved.
\[
\]
As someone can see these quantities are behave as (RR), the (RC) and (RGG) continued fractions identities.
For example when $q=e^{-\pi\sqrt{r}}$, $r$ positive rational, they are algebraic numbers and satisfy modular equations. Their derivatives also are all obey the same nome.\\         Let now
\begin{equation}
\left(a;q\right)_k:=\prod^{k-1}_{n=0}(1-aq^n)
\end{equation}The Rogers Ramanujan continued fraction is
\begin{equation}
R(q):=\frac{q^{1/5}}{1+}\frac{q}{1+}\frac{q^2}{1+}\frac{q^3}{1+}\cdots  
\end{equation}
which satisfies the famous Rogers-Ramanujan identity:
\begin{equation}
R^{*}(q):=q^{-1/5}R(q)=\frac{(q;q^5)_{\infty}(q^4;q^5)_{\infty}}{(q^2;q^5)_{\infty}(q^3;q^5)_{\infty}}=\prod^{\infty}_{n=1}(1-q^n)^{X_2(n)}
\end{equation}
where $X_2(n)$ is the Legendre symbol $\left(\frac{n}{5}\right)$.\\Also hold
\begin{equation}
R(e^{-x})=e^{-x/5}\frac{\vartheta_4(3ix/4,e^{-5x/2})}{\vartheta_4(ix/4,e^{-5x/2})},  x>0
\end{equation} 
Where $\vartheta_4(a,q)$ is the 4th kind Elliptic Theta function (see [9]).\\
The concept of formulation (1) is described below. We first begin with the rewriting of (4) into the form (see [14])
\begin{equation}
R(e^{-x})=\exp\left(-x/5-\sum^{\infty}_{n=1}\frac{1}{n}\frac{e^{4nx}-e^{3nx}-e^{2nx}+e^{nx}}{e^{5nx}-1}\right), x>0
\end{equation}
The Ramanujan-Gollnitz-Gordon continued fraction is 
\begin{equation}
H(q)=\frac{q^{1/2}}{(1+q)+}\frac{q^2}{(1+q^3)+}\frac{q^4}{(1+q^5)+}\frac{q^6}{(1+q^7)+}\cdots
\end{equation} 
Also for this continued fraction holds  
\begin{equation}
H(e^{-x})=\exp\left(-x/2-\sum^{\infty}_{n=1}\frac{1}{n}\frac{e^{7nx}-e^{5nx}-e^{3nx}+e^{nx}}{e^{8nx}-1}\right), x>0
\end{equation}
\begin{equation}
H(e^{-x})=e^{-x/2}\frac{\vartheta_4(3ix/2,e^{-4x})}{\vartheta_4(ix/2,e^{-4x})},  x>0
\end{equation}
\[
\]
Is true that exists generalizations for these expansions, but there is no theory developed, especially for evaluations and modular equations.

\section{Theorems on Rogers Ramanujan Quantities}

\textbf{Definition 1.} In general if $q=e^{-\pi \sqrt{r}}$ where $a,p,r>0$ we denote 'Agile' the quantity 
\begin{equation}
[a,p;q]=(q^{p-a};q^p)_{\infty}(q^a;q^p)_{\infty}
\end{equation}
\textbf{Definition 2.} We call  
\begin{equation}
R(a,b,p;q):=q^{-(a-b)/2+(a^2-b^2)/(2p)}\frac{[a,p;q]}{[b,p;q]}
\end{equation}
'Ramanujan's Quantity' because many of Ramanujan's continued fractions can be put in this form.\\
Also
$$
R^{*}(a,b,p;q):=\frac{[a,p;q]}{[b,p;q]}
$$
\[
\] 
\textbf{Conjecture 1.} 
If $q=e^{-\pi \sqrt{r}}$, $a,b,p,r$ positive rationals then
\begin{equation}
q^{p/12-a/2+a^2/{(2p)}}[a,p;q]\stackrel{?}{=}Algebraic
\end{equation}
\textbf{Note.}
The mark ''?'' means that we have no proof.
\[
\]
\textbf{Lemma 1.}
\begin{equation}
\sum^{\infty}_{k=1}\frac{\cosh(2tk)}{k\sinh(\pi ak)}=\log(P_0)-\log(\vartheta_4(it,e^{-a\pi}))\textrm{ , where } \left|2t\right|<\left|\pi a\right|
\end{equation}
and $P_0=\prod^{\infty}_{n=1}(1-e^{-2n\pi a})$ and $\vartheta_4(u,q)=1+\sum^{\infty}_{n=1}(-1)^nq^{n^2}\cos(2nu)$.\\
\textbf{Proof.}\\ From ([2] pg.170 relation (13-2-12)) and the definition of theta functions we have
\begin{equation}
\vartheta_4(z,q)=\prod^{\infty}_{n=0}(1-q^{2n+2})(1-q^{2n-1}e^{2iz})(1-q^{2n-1}e^{-2iz})
\end{equation}
By taking the logarithm of both sides and expanding the logarithm of the individual terms in a power series it is simple to show (13) from (14), where $q=e^{-\pi a}$, $a$ positive real.
\[
\] 
\textbf{Theorem 1.} 
If $a,b,p,r$ are positive rationals, then
\begin{equation}
R(a,b,p;q):=q^{-(a-b)/2+(a^2-b^2)/(2p)}R^{*}(a,b,p;q)=Algebraic
\end{equation} 
\textbf{Proof.}\\
Eq.(15) follows easy from the Conjecture 1 and the Definitions 1,2.
\[
\]
One example is the Rogers-Ramanujan continued fraction
\begin{equation}
q^{1/5}R^*(1,2,5;q)=R^{*}(q)q^{1/5}=R(q)
\end{equation}
\[
\]
\textbf{Theorem 2.} For all positive reals $a,b,p,x$ 
\begin{equation} R(a,b,p;e^{-x})=\exp\left(-x\frac{a^2-b^2}{2p}+x\frac{a-b}{2}\right)\frac{\vartheta_4\left((p-2a)ix/4,e^{-px/2}\right)}{\vartheta_4\left((p-2b)ix/4,e^{-px/2}\right)}=
\end{equation}
\begin{equation}
=\exp\left[-x\left(\frac{a^2-b^2}{2p}-\frac{a-b}{2}\right)-\sum^{\infty}_{n=1}\frac{1}{n}\frac{e^{anx}+e^{(p-a)nx}-e^{(p-b)nx}-e^{bnx}}{e^{pnx}-1}\right]
\end{equation}
\textbf{Proof.}\\ 
From Definitions 1, 2 and the relations (13), (14) we can rewrite $R$ in the form 
$$
R(a,b,p;e^{-x})=\exp\left(-x\frac{a^2-b^2}{2p}+x\frac{a-b}{2}\right)\frac{\exp\left(\sum^{\infty}_{n=1}\frac{\cosh(nx(p-2b)/2)}{n \sinh(pnx/2)}\right)}{\exp\left(\sum^{\infty}_{n=1}\frac{\cosh((p-2a)nx/2)}{n \sinh(pnx/2)}\right)}
$$
from which as one can see (17) and (18) follow.  
\[
\]
For the continued fraction (7) we give some evaluations with the command 'Recognize' of Mathematica:
$$
H\left(e^{-\pi}\right)\stackrel{?}{=}\sqrt{4-2\sqrt{2}}-1-\sqrt{2}
$$
$$
H\left(e^{-\pi\sqrt{2}}\right)\stackrel{?}{=}(1-8t-12t^2-8t^3+38t^4+8t^5-12t^6+8t^7+t^8)_3
$$
Where $(P(x))_n$ is the nth root of the equation $P(x)=0$ taken with Mathematica program. 
\[
\]
\textbf{Theorem 3.} (The Rogers Ramanujan Identity of the Quantities)\\ If $a,b,p$ are positive integers and $p-a\neq p-b$, $\left|q\right|<1$, then
\begin{equation}
R^{*}(a,b,p;q)=\prod^{\infty}_{n=1}(1-q^n)^{X(n)}
\end{equation}
where
\begin{equation}
X(n)= \left\{\begin{array}{cc}
                                    1, \mbox{  }  n\equiv(p-a)modp \\
                                   -1, \mbox{  }  n\equiv(p-b)modp \\
                                    1, \mbox{  }  n\equiv a mod p \\
                                   -1, \mbox{  }  n\equiv bmodp \\
                                    0, \mbox{  }  else 
\end{array}\right\}
\end{equation}
\textbf{Proof.}\\ 
Use Theorem 2. Take the logarithms and expand the product (19). The proof is easy.
\[
\]
Hence if we set 
\begin{equation}
M(q)=M(a,b,p;q):=q\frac{dR(a,b,p;q)}{dq}\frac{1}{R(a,b,p;q)}
\end{equation}
Then 
\begin{equation}
q\frac{dR(a,b,p;q)}{dq}\frac{1}{R(a,b,p;q)}=Q-\sum^{\infty}_{n=1}q^n\tau(n) , 
\end{equation} 
where 
\begin{equation}
Q=-(a-b)/2+(a^2-b^2)/(2p)  
\end{equation}
\begin{equation}
\tau(n)=\tau(a,b,p,n):=\sum_{d|n}X(d)d .
\end{equation}
\[
\]
\textbf{Theorem 4.} Let $|q|<1$, then
\begin{equation}
\log(R^{*}(a,b;p;q))=-\sum^{\infty}_{n=1}\frac{q^n}{n}\sum_{d|n}{X(d)d}
\end{equation}
\textbf{Proof.}\\
Follows from Theorem 3.
\[
\]
\textbf{Theorem 5.} 
For every positive integers $a$, $b$, $p$ with $a<p$ and $b<p$, with $p$ prime we have
\begin{equation}
\tau(n)=\tau(pn)
\end{equation}
where $n=1,2,3,\ldots$\\
\textbf{Proof.}\\
For $n=1$, then $\tau(1)=X(1)$ and $\tau(p)=X(1)1+X(p)p=X(1)$, (from (20)).    
Let for $n=n_0$ is $\tau(n_0)=\sum_{d|n_0}X(d)d$.
For $n_1=n_0p$, $p$-prime: $\tau(n_0p)=\sum_{d|(n_0p)}X(d)d$. 
Let $(p,n_0)=1$, then $\tau(n_0)=\tau(n_0p)$ or equivalently $\sum_{d|n_0}X(pd)pd=0$, which is true since $p>a,b$. If was $(p,n_0)=s>1$, then $s=p$ and exist $k_0$, $m$: $n_0=mp^{k_0}$, $(m,p)=1$. Again all the multiplies of $p$ are canceled and we lead to $\sum_{d|mp^{k_{0}}}X(d)d-\tau(m)=0$. Hence we get the proof.
\[
\]
\textbf{Proposition.}(See [7] pg. 24)\\ Suppose that $a,b$ and $q$ are complex numbers with $\left|a b \right|<1$ and $\left|q\right|<1$ or that $a=b^{2m+1}$ for some integer $m$. Then  
$$
P(a,b,q):=\frac{(a^2 q^3;q^4)_{\infty} (b^2q^3;q^4)_{\infty}}{(a^2q;q^4)_{\infty}(b^2q;q^4)_{\infty}}=$$
\begin{equation}
=\frac{1}{(1-ab)+}\frac{(a-bq)(b-aq)}{(1-ab)(q^2+1)+}\frac{(a-bq^3)(b-aq^3)}{(1-ab)(q^4+1)+}\frac{(a-bq^5)(b-aq^5)}{(1-ab)(q^6+1)+\ldots}
\end{equation}
\[
\]
\textbf{Theorem 6.} 
If $a=2A+3p/4$, $b=2B+p/4$ and $p=4(A+B)$, $|q|<1$
\begin{equation}
R^*(a,b,p;q)=(1-q^{B-A})P(q^A,q^B,q^{A+B}) 
\end{equation}
\textbf{Proof.}\\
One can see that
\begin{equation}
P(q^A,q^B,q^{A+B})=\frac{(q^a;q^p)_{\infty}(q^{2p-a};q^p)_{\infty}}{[b,p;q]} 
\end{equation}
where $a=2A+3p/4$, $b=2B+p/4$ and $p=4(A+B)$.\\
Define 
\begin{equation}
{}_2\phi_1[a,b;c;q,z]:=\sum^{\infty}_{n=0} \frac{(a;q)_n (b;q)_n}{(c;q)_n}\frac{z^n}{(q;q)_n}
\end{equation}
and
\begin{equation}
\psi(a,q,z):=\sum^{\infty}_{n=0}\frac{(a;q)_n}{(q,q)_n}z^n={}_2\phi_1[a,0,0,q,z]
\end{equation}
Then
\begin{equation}
\psi(q^p,q^p,q^{p-a})R^*(a,b,p;q)=P[q^A,q^B,q^{A+B}] 
\end{equation}
The proof of (28) follows easily from (32) and the $q$-binomial theorem (see [7]): 
\[
\psi(a,q,z)=\prod^{\infty}_{n=0}\frac{1-azq^n}{1-zq^n}
\]
\textbf{Note.}
Relation (28) is an expansion of a Ramanujan Quantity in continued fraction.\\
Also in view of [7] page 14 Entry 4, we have 
$$
R^{*}(a,b,p;q)\sum^{\infty}_{n=0}\frac{(q^{a-b};q^p)_n(q^{b+a-p};q^p)_n}{(q^p;q^p)_n(q^a;q^p)_n}q^{(p-a)n}=1 . 
$$
\textbf{Theorem 7.} 
If $a,b,p\in\bf N\rm$, with $a<b<p$, $p-$prime and $\zeta^m_p=e^{-\pi i m/p}$, then
\begin{equation} 
R(a,b,p;q^p)=\prod^{p-1}_{m=0}R\left(a,b,p;(-1)^m\zeta^m_pq\right) 
\end{equation}
\textbf{Proof.}\\
From Theorem 5 holds $\tau(n)=\tau(np)$, hence 
\begin{equation}
\log(R(a,b,p;q))=Q\log(q)-\sum^{\infty}_{n=1}\frac{q^{n }}{n}\tau(n)=
Q\log(q)-p\sum^{\infty}_{n=1}\frac{q^{n}}{n p}\tau(np)
\end{equation}
Also holds the general identity:\\ 
If
$$
y(x)=\sum^{\infty}_{n=0}a_nx^n$$ then
\begin{equation}  \sum^{\infty}_{n=0}a_{np}x^{np}=p^{-1}\sum^{p-1}_{m=0}y\left((-1)^m\zeta^m_{p}x\right)
\end{equation}  
From (34) and (35) we get the proof.
\[
\]
\textbf{Note.} 
Another form of Theorem 7 is: 
If $a,b\in\bf N\rm$ and $p$-positive prime: $a, b<p$, then 
\begin{equation}
R(a,b,p;q^p)=\prod^{p-1}_{m=0}R\left(a,b,p;e^{2\pi i m/p}q\right)
\end{equation}
\[
\]
\textbf{Conjecture.}\\
i) If $a$, $b$ are odd positive integers and $p$-even, then 
\begin{equation}
\tau(a,b,p,2n)\stackrel{?}{=}\tau(a,b,p,n)
\end{equation}
ii) If $a\neq 0(mod3)$, $b\neq 0(mod3)$ and $p$ positive integer such that $(p,3)>1$, then 
\begin{equation}
\tau(a,b,p,3n)\stackrel{?}{=}\tau(a,b,p,n)
\end{equation} 
\ldots
\[
\]
In general if $p_0$ is prime and if $a\neq 0(modp_0)$, $b\neq 0(modp_0)$ and $p$ is positive integer such that $(p,p_0)>1$ then 
\begin{equation}
\tau(a,b,p,p_0 n)\stackrel{?}{=}\tau(a,b,p,n)
\end{equation} 
\[
\]
\textbf{Theorem 8.} 
If $a,b,p$ as in Conjecture-(i) and $R(q)=R(a,b,p;q)$, then
\begin{equation}
2M(q^2)\stackrel{?}{=}M(q)+M(-q)
\end{equation}
\begin{equation}
R(q^2)\stackrel{?}{=}R(q)\left|R(-q)\right|
\end{equation}
\textbf{Proof.}\\
From $\tau(2n)=\tau(n)$ we have, expanding $M(q)$ into series: 
$$
M(-q)+M(q)=2Q-2\sum^{\infty}_{n=1}\tau(2n)q^{2n}=2\left(Q-\sum^{\infty}_{n=1}\tau(n)(q^2)^n\right)=2M(q^2)
$$ 

\section{The Programs}

In this section we will give some results with no proof. Also we present the program which we use to found these $\tau$ relations. 
\[
\] 
\centerline{\textbf{ The program of $\tau$ with Mathematica}}
\[
\]
\textbf{Routine 1.}\\
Suppose $a=1$, $b=4$, $p=17$. Set
$${\rm Q1[n,k]:=If[Mod[n,k]==0,1,0]}$$ 
$${\rm t[n]:=\sum^{n}_{k=1}X[1,4,17,k]kQ1[n,k]}$$  
then use the 'Solve' routine as
$$
{\rm Solve[Table[\sum^{17}_{j=1}c[j ]t[jn]==0,\left\{n,1,289\right\}],Table[c[r],\left\{r,1,17\right\}]]} 
$$
We get
$$
\rm{\{\{c[1]\to -4 c[8]-4 c[16]-c[17],c[2]\to 4 c[8],c[3]\to 0,c[4]\to -c[8]+3 c[16],}$$
$$\rm{,c[5]\to 0,c[6]\to 0,c[7]\to 0,c[9]\to 0,c[10]\to 0,c[11]\to 0,c[12]\to 0,}$$
$$\rm{,c[13]\to 0,c[14]\to 0,c[15]\to 0\}\}}
$$
We change in the above output equalities the symbol $\to$ into $=$ 
$${\rm Clear[t]}$$
$$\rm{Coefficient[\sum^{17}_{j=1}c[j]t[j n],c[16]]}$$
$$\rm{-4 t[n]+3 t[4 n]+t[16 n]}$$ 
$$\rm{Coefficient[\sum^{17}_{j=1}c[j]t[j n],c[8]]}$$
$$\rm{-4 t[n]+4 t[2 n]-t[4 n]+t[8 n]}$$ 
\[
\]
With the above Program one can find for every Ramanujan Quantity relations between the $\tau$'s. We present such  relations in the next two Propositions:
\[
\]
\textbf{Proposition 1.} 
If $a=1$, $b=4$, $p=17$ we have
\begin{equation}
-4\tau(n)+3\tau(4n)+\tau(16n)=0
\end{equation}
and
\begin{equation}
-4\tau(n)+4\tau(2n)-\tau(4n)+\tau(8n)=0
\end{equation} 
\[
\]
\textbf{Proposition 2.} 
If $a=1$, $b=5$, $p=26$ we have
\begin{equation}
-\tau(n)+\tau(3n)-\tau(5n)+\tau(15n)=0
\end{equation}
\begin{equation}
-\frac{26}{77}\tau(n)-\frac{17}{7}\tau(3n)+\frac{17}{7}\tau(7n)-\frac{51}{77}\tau(11n)+\tau(17n)=0
\end{equation} 
\begin{equation}
-\frac{134}{77}\tau(n)+\frac{19}{7} \tau(3 n)-\frac{19}{7} \tau(7 n)+\frac{57}{77} \tau(11 n)+\tau(19 n)=0
\end{equation}
\begin{equation}
-\frac{34}{11}\tau(n)+\frac{23}{11} \tau(11 n)+\tau(23 n)=0
\end{equation}
\begin{equation}
-5\tau(n)+4\tau(5n)+\tau(25n)=0
\end{equation}
\[
\]
\centerline{\textbf{ A simple program for finding possible modular equations}}\\
\centerline{\textbf{with Mathematica}}
\[
\]
\textbf{Routine 2.}
$$
x=\textrm{Series}[\textrm{R}[a,b,p,q^{\alpha}],\{q,0,A\}]
$$
$$
y=\textrm{Series}[\textrm{R}[a,b,p,q^{\beta}],\{q,0,A\}]
$$
$s=\ldots$
$$
t=\textrm{Table}\left[\textrm{Coefficient}\left[\sum _{i=0}^s \sum _{j=0}^s c[i,j] x{}^{\wedge}i y{}^{\wedge}j,q{}^{\wedge}n\right]==0,\{n,1,s{}^{\wedge}2\}\right];
$$
$$
r=\textrm{Table}[c[i,j],\{i,0,s\},\{j,0,s\}];
$$
$$
\textrm{Clear}[u,v]
$$
$$
\textrm{r1}=\textrm{Table}[u{}^{\wedge}i v{}^{\wedge}j,\{i,0,s\},\{j,0,s\}];
$$
$$
m=\textrm{Normal}[\textrm{Extract}[\textrm{CoefficientArrays}[t\textrm{//}\textrm{Flatten},r\textrm{//}\textrm{Flatten}],2]];
$$
$$
\textrm{m0}=\textrm{Normal}[m];
$$
$$
\textrm{Take}[\textrm{NullSpace}[\textrm{m0}],1] . \textrm{Flatten}[r]
$$
$$
\textrm{Take}[\textrm{NullSpace}[\textrm{m0}],1] . \textrm{Flatten}[\textrm{r1}]
$$

\section{The first Order Derivatives of Ramanujan's Quantities}

Observe that if
\begin{equation}
R_1(q)=\frac{q^{1/5}}{1+}\frac{q}{1+}\frac{q^2}{1+}\frac{q^3}{1+}\cdots=q^{1/5}\frac{(q;q^5)_{\infty}(q^4;q^5)_{\infty}}{(q^2;q^5)_{\infty}(q^3;q^5)_{\infty}}
\end{equation}    
\begin{equation}
R_2(q)=\frac{q^{1/3}}{1+}\frac{q+q^2}{1+}\frac{q^2+q^4}{1+}\frac{q^3+q^6}{1+}\cdots=q^{1/3}\frac{(q;q^6)_{\infty}(q^5;q^6)_{\infty}}{(q^3;q^6)_\infty^2}
\end{equation}    
\begin{equation}
R_3(q)=\frac{q^{1/2}}{(1+q)+}\frac{q^2}{(1+q^3)+}\frac{q^4}{(1+q^5)+}\cdots=q^{1/2}\frac{(q;q^8)_{\infty}(q^7;q^8)_{\infty}}{(q^3;q^8)_\infty (q^5;q^8)_{\infty}}
\end{equation}
\[
\]    
$R_1(q)=R(q)=R(1,3,5;q)$, $R_2(q)=V(q)=R(1,3,6;q)$, $R_3(q)=H(q)=R(1,3,8;q)$ are respectively the Rogers-Ramanujan, Ramanujan's Cubic and Ramanujan-Gollnitz-Gordon continued fraction, these have derivatives
\begin{equation} 
R_{1,2,3}'(q)\frac{q\pi^2}{K(k_r)^2}\stackrel{?}{=}Algebraic
\end{equation}
whenever $q=e^{-\pi \sqrt{r}}$ and $r$ is a positive rational.
\[
\]
\textbf{Conjecture 2.} If $a,b,p,r$ are positive rationals with $a,b<p$, then
\begin{equation}
\frac{d}{dq}R(a,b,p;q)\stackrel{?}{=}\frac{K(k_r)^2}{q \pi^2}Algebraic
\end{equation}
\begin{equation}
\frac{d}{dq}\left(q^{p/12-a/2+a^2/(2p)}[a,p;q]\right)\stackrel{?}{=}\frac{K(k_r)^2}{q \pi^2}Algebraic
\end{equation}
\[
\]
Let now $k=k_r$ be the Elliptic singular modulus (see [10],[12],[15]). 
In [12] we have prove the following relation (see also [6] pg. 87]), for real $r>0$:
\begin{equation}
\frac{dr}{dk}=\frac{\pi\sqrt{r}}{K(k)^2k(1-k^2)}
\end{equation}
Hence
\begin{equation}
\frac{dq}{dk}=\frac{-q\pi^2}{2k(1-k^2)K(k)^2}
\end{equation}
This observation along with Conjecture 2 lead us to the concluding remark
$$
\frac{dR(a,b,p;q)}{dk}=$$
$$=\frac{dR(a,b,p;q)}{dq}\frac{dq}{dk}=\frac{K(k)^2}{q\pi^2}\frac{-q \pi^2}{2K(k)^22k(1-k^2)}\cdot Algebraic$$
Hence
\[
\]
\textbf{Theorem 9.} 
When $q=e^{-\pi\sqrt{r}}$, $a,b,p,r$ positive rationals, then 
\begin{equation}
\frac{dR(a,b,p;q)}{dk}\stackrel{?}{=}Algebraic
\end{equation}
\[
\]
\textbf{Theorem 10.} 
If $q=e^{-\pi\sqrt{r}}$, then
\begin{equation}
\frac{dH(q)}{dq}=\frac{-q\pi^2}{2k(1-k^2)K(k)^2}\frac{\sqrt{1-k'}}{k'(k\sqrt{2}+2\sqrt{1-k'})} 
\end{equation}
\textbf{Proof.}\\
In [10], we have proved that
\begin{equation}
H(q)=-t+\sqrt{t^2+1}\textrm{ , }t=\frac{k_r}{(1-k'_r)}
\end{equation}
which gives
\begin{equation}
\frac{dH(q)}{dk}=\frac{\sqrt{1-k'}}{k'(k\sqrt{2}+2\sqrt{1-k'})} 
\end{equation} 
using now (56) we get the result.
\[
\]
\textbf{Note.} In [12] we have derive the first order derivative for the Cubic Continued fraction:\\
Let $q=e^{-\pi\sqrt{r}}$, $r>0$ then
\begin{equation}
V'(q)=\frac{dV(q)}{dq}=\frac{4K^2(k_r)k'^2_r(V(q)+V^4(q))}{3q\pi^2\sqrt{r}\sqrt{1-8V^3(q)}}
\end{equation}  
Hence we have prove that (50) and (51) obey (52). In [10] we have given a formula for the (RR) first derivative involving equations that can not solve in radicals (higher than 4). Also in [13] we have given a formula for (RR) but contains the function $k^{(-1)}(x)$, which is the inverse function of $k_{r}$.     
\[
\] 
\textbf{Examples.}
\begin{equation}
\left(\frac{d}{dq}R(1,2,4;q)\right)_{q=e^{-\pi}}\stackrel{?}{=}\frac{e^{\pi}\Gamma(1/4)^4}{64\cdot2^{5/8}\pi^3}
\end{equation}
\begin{equation}
\left(\frac{d}{dq}R(1,2,5;q)\right)_{q=e^{-\pi}}\stackrel{?}{=}\frac{e^{\pi}\Gamma(1/4)^4}{16\pi^3}p_1
\end{equation}
where
\[
p_1=(16-240t^2+800t^3-2900t^4-6000t^5-6500t^6+17500t^7+625t^8)_3
\]
\begin{equation}
\left(\frac{d}{dq}R(1,3,8;q)\right)_{q=e^{-\pi}}\stackrel{?}{=}\left(2+\sqrt{2}-\sqrt{5-\frac{7}{2}\sqrt{2}}\right)\frac{64 e^{\pi } \pi }{\Gamma\left(-\frac{1}{4}\right)^4}
\end{equation}
\begin{equation}
\left(\frac{d}{dq}R(1,3,8;q)\right)_{q=e^{-2\pi}}\stackrel{?}{=}\frac{(6+4\sqrt{2})e^{2 \pi}\Gamma(5/4)^4}{\pi^3}p_2
\end{equation}
\[
p_2=(16384-1720320 t^2-6684672 t^3+143104 t^4-18432 t^5-1664 t^6+t^8)_3
\]
\[
\] 

\section{Modular equations and Ramanujan Quantities}

With the help of Theorem 2 we can evaluate $R(a,b,p;q)$ in series of $q^{Q}$: $$R(a,b,p;q)=\sum^{M}_{n=0}c_nq^{nQ} , \eqno{(a)}$$ where $M$ is positive integer and 
$$Q=\frac{a^2-b^2}{2p}-\frac{a-b}{2}$$
Setting as in [11]: $$R_S=\sum_{0\leq i+j\leq d}a_{ij}u^iv^j $$ and $d$ is suitable large positive integer, we try to solve $R_S=0$, where $$u=R(a,b,p;q) , v=R(a,b,p;q^{\nu})  $$ are given from (a) and $\nu$ positive integer. Evaluating the $a_{ij}$, we obtain  modular equations for $R(a,b,p;q)$.
\[
\]    
1) We present some modular equations for the Ramanujan Quantity $R(1,2,4;q)$, which:
$$
R(1,2,4;q)=1+\frac{(1+q)}{1+}\frac{q^2}{1+}\frac{(q+q^3)}{1+}\frac{q^4}{1+}\ldots=\frac{(q^2;q^4)^2_{\infty}}{(q;q^4)_{\infty}(q^3;q^4)_{\infty}}
$$
\[
\]
a) If $u=R(1,2,4;q)$ and $v=R(1,2,4;q^2)$, then 
\begin{equation} 
u^4-v^2+4u^4v^4\stackrel{?}{=}0
\end{equation}
b) If $u=R(1,2,4;q)$ and $v=R(1,2,4;q^3)$, then
\begin{equation}
u^4-uv+4u^3v^3-v^4\stackrel{?}{=0}
\end{equation}
c) If $u=R(1,2,4;q)$ and $v=R(1,2,4;q^5)$, then
\begin{equation}
u^6-uv+5u^4v^2-5u^2v^4+16 u^5 v^5-v^6\stackrel{?}{=0}
\end{equation}
d) If $u=R(1,2,4;q)$ and $v=R(1,2,4,q^7)$, then
\begin{equation}
u^8-u v+7 u^2 v^2-28 u^3 v^3+70 u^4 v^4-112 u^5 v^5+112 u^6 v^6-64 u^7 v^7+v^8\stackrel{?}{=0}
\end{equation}
\[
\] 
2) For the Ramanujan Quantity $R(1,2,6;q)$ we have
\[
\]
a) If $u=R(1,2,6;q)$ and $v=R(1,2,6;q^2)$, then
\begin{equation}   
u^4-v^2+3u^4 v^2+v^4\stackrel{?}{=0}
\end{equation}
$\ldots$\\
One can find with the help of Mathematica many relations such above 
\[
\]
\textbf{The 5-degree modular equation of Ramanujan's Cubic continued fraction:}\\ 
If $u=R(1,3,6;q)$ and $v=R(1,3,6;q^5)$, then
\begin{equation}
u^6-uv+5u^4v+5u^2v^2-10u^5v^2-20u^3v^3+5uv^4+20u^4v^4-10u^2 v^5-16u^5v^5+v^6=0
\end{equation}
\textbf{The 7-degree modular equation of Ramanujan's Cubic continued fraction:}\\ 
If $u=R(1,3,6;q)$ and $v=R(1,3,6;q^7)$, then
$$u^8-uv+7u^4v+28u^6v^2-56u^5v^3+7uv^4+21u^4v^4-56u^7v^4-56u^3 v^5+$$
\begin{equation}
+28u^2v^6-56u^4v^7-64u^7v^7+v^8=0
\end{equation}
For proofs see at the end of page 15.
\[
\]
If $a>b$ then from the definition of the Ramanujan Quantity (RQ) we have 
\begin{equation}
R(a,b,p;q)=\frac{1}{R(b,a,p;q)}
\end{equation} 
Suppose that $a=\frac{a_1}{a_2}$, $b=\frac{b_1}{b_2}$, $p=\frac{p_1}{p_2}$, and $u(q)=R(a,b,p;q)$, then $$u(q^{\frac{1}{a_2b_2p_2}})=R(a_1b_2p_2,b_1a_2p_2,p_1a_2b_2;q)=w(q), $$ if $a_1b_2p_2<b_1a_2p_2$, (otherwise we use (73)).\\ But $w_1:=w(q)$, $w_{\nu}:=w(q^{\nu})$ are related by a modular equation $f(w_1,w_{\nu})=0$, or $f(w(q^{a_2b_2p_2}),w(q^{\nu\cdot a_2b_2p_2}))=0$. Hence
\[
\]
\textbf{Theorem 11.} 
When $a=\frac{a_1}{a_2}$, $b=\frac{b_1}{b_2}$, $p=\frac{p_1}{p_2}$, $a_1$, $a_2$, $b_1$, $b_2$, $p_1$, $p_2\in\bf N\rm$ and $a_1b_2p_2<b_1a_2p_2$ then the modular equation which relates $u_1:=R(a,b,p;q)$ and $u_{\nu}:=R(a,b,p;q^{\nu})$, $\nu\in \bf N\rm $ is that of 
\begin{equation}
w(q)=R(a_1b_2p_2,b_1a_2p_2,p_1a_2b_2;q) \textrm{ and } w(q^{\nu}).
\end{equation}
\[
\]
\textbf{Example.}\\
The modular equation between $z_1=z(q)=R\left(1,\frac{1}{2},2;q\right)$ and $z_2=z(q^2)$ is 
\begin{equation}
4+z_2^4-z_1^4z_2^2\stackrel{?}{=}0
\end{equation}   
\textbf{Proof.}\\
We have $$z_1=z(q)=R\left(1,\frac{1}{2},2;q\right)=R(2,1,4;q^{1/2})=\frac{1}{R(1,2,4;q^{1/2})}=\frac{1}{u(q^{1/2})}$$
But $z_2=z(q^2)=\frac{1}{u(q^{2/2})}$ using (66) we have  
$$u(q^{1/2})^4-u(q)^2+4u(q^{1/2})^4u(q)^4=0 , $$
from which (75) follows.
\[
\] 
From Theorem 4 differentiating (25) and using Conjecture 2 we have that
\begin{equation}
q\frac{dR(a,b,p;q)}{dq}\frac{1}{R(a,b,p;q)}=Q-\sum^{\infty}_{n=1}q^n\sum_{d|n}X(d)d
\end{equation} 
along with (see [14]): 
$$
f(-q)^4=2^{4/3}\pi^{-2}q^{-1/6}(k_r)^{1/3}(k'_r)^{4/3}K(k_r)^2\eqno{:(b)}
$$
\begin{equation}
N(q)=q^{-1/6}f(-q)^{-4}\left(Q-\sum^{\infty}_{n=1}q^n\sum_{d|n}X(d)d\right)=Algebraic 
\end{equation}
This is a resulting formula for the first derivative:
\begin{equation}
q^{5/6}f(-q)^{-4}\frac{R'(q)}{R(q)}=N(q)  
\end{equation}
The function $N(q)$ take algebraic values when $q=e^{-\pi\sqrt{r}}$, $r$ positive rational and in the case of (RR) satisfies modular equations. With the same method as in $R(a,b,p;q)$ which we use in the beginning of the paragraph 4 we have:
\[
\]
\textbf{The 2-degree Modular equation for the first derivative of RR continued fraction}
\[
\] 
For $a=1$, $b=2$, $p=5$, we have the case of (RR) and\\ 
a) If $u=N(q)$ and $v=N(q^2)$ then     
\begin{equation}
5u^6-u^2v^2-125u^4v^4+5v^6=0
\end{equation}
\[
\]
\textbf{The 3-degree Modular equation for the first derivative of RR continued fraction}
\[
\] 
b) If $u=N(q)$ and $v=N(q^3)$ we have
\begin{equation}
125u^{12}+u^3v^3+1125u^9v^3+1125u^3v^9+1953125u^9v^9-125v^{12}=0
\end{equation}
\[
\]
Proofs for (79) and (80) can given if we use the identity found in [16]:
$$
R'(q)=5^{-1}q^{-5/6}f(-q)^4R(q)\left[R(q)^{-5}-11-R(q)^5\right]^{-1/6}
$$
and Ramanujan's duplication and triplication formulas of $R(q)$:
$$
\frac{R(q^2)-R(q)^2}{R(q^2)+R(q)^2}=R(q)R(q^2)^2,
$$
$$
\left[R(q^3)-R(q)^3\right]\left[1+R(q)R(q^3)^3\right]=3R(q)^2R(q^3)^2.
$$
The hard thing is that have very tedious algebraic calculations. In view of [12] the same ideas hold and for the modular equations (71) and (72).
\[
\]
Suppose now that $q_0=e^{-\pi\sqrt{r_0}}$ and we know  $R^{(1)}(q_0)=\left(\frac{dR(q)}{dq}\right)_{q=q_0}$, then from equations (76),(77),(78),(79) and (80) we can evaluate in radicals, any high order values of the first derivative of the (RR) in which $r=4^n9^mr_0$, for $n$, $m$ integers.
\[
\] 
\textbf{Note.}
If $K(k_r)=K[r]$ then holds: 
$$K[4r]=\frac{1+k'_r}{2}K[r]$$
and
$$K[9r]=m_3(r)K[r]$$
where $m_3(r)$ is solution of       
$$27m_3(r)^4-18m_3(r)^2-8(1-2k_r^2)m_3(r)-1=0$$
The formulas for evaluation of $k_{4r}$ and $k_{9r}$ are in [7] and [15].
\[
\]
\textbf{Proposition.}
If $R(q)$ is the Rogers Ramanujan continued fraction and $$M(q)=\frac{R'(q)}{R(q)}$$ then 
$$M(q^2)=\frac{\sqrt{q^{5/2}f(-q^2)^8M(q)\left(125q^{5/2}M(q)-\sqrt{5}\sqrt{8f(-q)^{24}+3125q^5M(q)}\right)}}{2\sqrt{5}q^{5/2}f(-q)^8}$$
where
$$
f(-q^2)=4\frac{(2m_2-1)^4}{m_2(m_2-1)}qf(-q)^{24}\textrm{ , where }m_2=\frac{1+k'_r}{2}
$$
\textbf{Proof.}\\
The proof follows from (79), the above Note and $(b)$ of page 14.

\section{Application in 'almost all' continued fractions}

In this section we give examples and evaluations of a certain class of continued fractions with values $a$, $b$, $p$ positive rationals under Theorem 6 holds.
\[
\]
\textbf{I.} The case of $A=1$, $B=2$,
then $a=11$, $b=7$, $p=12$ and  $$S(q)=S_{1,2}(q)=R(11,7,12;q)=\frac{1}{R(7,11,12;q)}$$ and from Theorem 5 we get 
$$
S_{1,2}(q)=
$$
$$
=q\frac{1-q}{1-q^3+}\frac{q^3 (1-q^2)(1-q^4)}{(1-q^3)(1+q^6)+}\frac{q^3(1-q^8)(1-q^{10})}{(1-q^3)(1+q^{12})+}\frac{q^3(1-q^{14})(1-q^{16})}{(1-q^3)(1+q^{18})+...}
$$
\begin{equation}
=\left(q\frac{\prod^{\infty}_{n=0}(1-q^{12n+11})(1-q^{12n+1})}{\prod^{\infty}_{n=0}(1-q^{12n+7})(1-q^{12n+5})}\right)^{-1}
\end{equation}
Hence with the above methods we find that continued fraction $S_{1,2}(q)$ obeys the following modular equations:\\  
1) If we set $u=S_{1,2}(q)$ and $v=S_{1,2}(q^2)$, then 
\begin{equation}
-u^2+v-2 u v+u^2 v-v^2\stackrel{?}{=0}
\end{equation}
2) If $u=S_{1,2}(q)$ and $v=S_{1,2}(q^3)$, then
\begin{equation}
u^3-v+3uv-u^3v+v^2-3u^2v^2+u^3v^2-v^3\stackrel{?}{=0}
\end{equation}
3) If $u=S_{1,2}(q)$ and $v=S_{1,2}(q^5)$, then
$$
-u^5+v-5 u v+5 u^2 v+5 u^5 v-10 u^3 v^2-5 u^5 v^2+10 u^2 v^3+10 u^4 v^3-5 u v^4-$$
\begin{equation}
-10 u^3 v^4+5 u v^5+5 u^4 v^5-5 u^5 v^5+u^6 v^5-u v^6\stackrel{?}{=0}
\end{equation}
4) If $u=S_{1,2}(q)$ and $v=S_{1,2}(q^7)$, then
$$     
-u^7+v-7uv+14u^2v-7u^3v+7u^5v-7u^6v+7u^7v+7u v^2-28u^2v^2+
$$
$$
+7u^3v^2-28u^5v^2+28u^6v^2-14u^7v^2-7uv^3+28u^2 v^3-7u^3v^3+
$$
$$
+35u^4v^3+7u^5v^3-7u^6v^3+7u^7v^3-35u^3v^4-35u^5 v^4+7uv^5-
$$
$$
-7u^2v^5+7u^3v^5+35u^4v^5-7u^5v^5+28u^6v^5-7u^7 v^5-14uv^6+
$$
$$
+28u^2v^6-28u^3v^6+7u^5v^6-28u^6v^6+7u^7v^6+7u v^7-7u^2v^7+
$$
\begin{equation}
+7u^3v^7-7u^5v^7+14u^6v^7-7u^7v^7+u^8v^7-u v^8\stackrel{?}{=0}
\end{equation}
\[
\]
\textbf{II.} In the case of $A=1$, $B=3$, then $a=14$, $b=10$, $p=16$ and $$S_{1,3}(q)=R(14,10,16;q)=\frac{1}{R(10,14,16;q)}$$ and from Theorem 5 we get 
$$
S_{1,3}(q)=
$$
$$=q\frac{1-q^2}{1-q^4+}\frac{q^4 (1-q^2)(1-q^6)}{(1-q^4)(1+q^8)+}\frac{q^4(1-q^{14})(1-q^{10})}{(1-q^4)(1+q^{12})+}\frac{q^4(1-q^{22})(1-q^{18})}{(1-q^4)(1+q^{20})+\ldots}=
$$
\begin{equation}
=\left(q\frac{\prod^{\infty}_{n=0}(1-q^{16n+10})(1-q^{16n+6})}{\prod^{\infty}_{n=0}(1-q^{16n+14})(1-q^{16n+2})}\right)^{-1}
\end{equation}
With the above methods we find that Continued fraction $S(q)$ obeys the following modular equations:\\  
1) If we set $u=S_{1,3}(q)$ and $v=S_{1,3}(q^2)$, then 
\begin{equation}
u^2-v+u^2 v+v^2\stackrel{?}{=0}
\end{equation}
2) If $u=S_{1,3}(q)$ and $v=S_{1,3}(q^3)$, then
\begin{equation}
u^3-v+3 u^2 v+3 u v^2-3 u^3 v^2-3 u^2 v^3+u^4 v^3-u v^4\stackrel{?}{=0}
\end{equation}
3) If $u=S_{1,3}(q)$ and $v=S_{1,3}(q^5)$, then
$$
u^5-v+5 u^2 v+10 u^3 v^2-5 u^5 v^2-10 u^2 v^3+10 u^4 v^3+5 u v^4-10 u^3 v^4-$$
\begin{equation}
-5 u^4 v^5+u^6 v^5-u v^6\stackrel{?}{=0}
\end{equation}
4) If $u=S_{1,3}(q)$ and $v=S_{1,3}(q^7)$, then
$$     
u^8-u v+7 u^3 v-7 u^5 v-7 u^7 v+28 u^6 v^2+7 u v^3-49 u^3 v^3-7 u^5 v^3-$$
$$
-7 u^7 v^3+70 u^4 v^4-7 u v^5-7 u^3 v^5-49 u^5 v^5+7 u^7 v^5+28 u^2 v^6-$$
\begin{equation}
-7 u v^7-7 u^3 v^7+7 u^5 v^7-u^7 v^7+v^8\stackrel{?}{=0}
\end{equation}
As a starting value we can get 
$$S_{1,3}(e^{-\pi/2})\stackrel{?}{=}-1-\sqrt{2}+\sqrt{2 \left(2+\sqrt{2}\right)}$$
from (87) we get 
$$S_{1,3}(e^{-\pi})=-3-2 \sqrt{2}+2 \sqrt{2+\sqrt{2}}+\sqrt{ 4+2\sqrt{2}}-$$
$$-\sqrt{30+22 \sqrt{2}-16 \sqrt{2+\sqrt{2}}-12 \sqrt{2 \left(2+\sqrt{2}\right)}}$$
$\ldots$ etc

A value for the derivative according to the Conjacture 2 is $$\left(\frac{dS_{1,3}(q)}{dq}\right)_{q=e^{-\pi/2}}\stackrel{?}{=} \frac{64 \left(4+2 \sqrt{2}-\sqrt{2 \left(10+7 \sqrt{2}\right)}\right) e^{\pi /2} \pi }{\Gamma\left(-\frac{1}{4}\right)^4}$$
\[
\]
\textbf{Note.} It is clear by this way that one can produce continued fractions such found and studied by Ramanujan.

\section{Approaches for the evaluation of Ramanujan Quantities and their derivatives}

In this section we will try to find Theoretical results that can be used to recover the Ramanujan Quantities and their derivatives. We will also drop the notation of $X$ which we use until now. Further
$$f(-q)=\prod^{\infty}_{n=1}(1-q^n)\eqno{(a)}$$
$$
L_1(q):=\sum^{\infty}_{n=1}\frac{nq^n}{1-q^n}\eqno{(b)}
$$
and $L(q)=1-24L_1(q)$.\\ The following elementary relations hold (see [7] pg. 120) 
$$
\sum^{\infty}_{n=1}\frac{(-1)^nnq^n}{1-q^n}=\frac{1}{24}(3-L(q)-4L(q^2))
$$  
$$
\sum^{\infty}_{n=1}\left(\frac{n}{4}\right)\frac{nq^n}{1-q^n}=\sum^{\infty}_{n=1}\left(\frac{n}{14}\right)\frac{nq^n}{1-q^n}=\frac{1}{24}(-1-L(q)+2L(q^2))
$$
$$
\left(\frac{n}{2}\right)=\left(\frac{n}{8}\right)\textrm{ , }\left(\frac{n}{4}\right)=\frac{1-(-1)^n}{2}
$$
$$
\sum^{\infty}_{n=1}\left(\frac{n}{3}\right)\frac{nq^n}{1-q^n}=\frac{1}{24}(2-L(q)-3L(q^3))+2\sum^{\infty}_{n=1}\frac{(3n+1)q^{3n+1}}{1-q^{3n+1}}
$$
We state a first result:
\[
\]
\textbf{Lemma 1.}
\begin{equation}
\frac{q}{R(a,b,p;q)}\frac{d}{dq}R(a,b,p;q)=M(q)=\sum^{\infty}_{n=1}\frac{X(n)nq^n}{1-q^{n}}
\end{equation}
\textbf{Proof.}\\
From the relation
$$
R^{*}(a,b,p;q)=\prod^{\infty}_{n=1}(1-q^n)^{X(a,b,p,n)} ,
$$
taking the logarithmic derivative, we get immediately the result. 
\[
\]
\textbf{Lemma 2.} 
It holds $\left(\frac{n}{g^2}\right)=1-X_g(n)$, where $g$ is prime and $X_g(n)$ the characteristic function of $0(\textrm{mod} g)$. Also hold 
\begin{equation}
\sum^{\infty}_{n=1}\left(\frac{n}{g^2}\right)\frac{nq^n}{1-q^n}=\frac{1}{24}(1-g-L(q)+gL(q^g)) 
\end{equation}
\textbf{Proof.}\\It follows easy from $\left(\frac{n}{g^2}\right)=1-X_g(n)$ and Lemma 1.
\[
\]
\textbf{Lemma 3.} 
If happens $X(a,b,p;n)=\left(\frac{n}{g^2}\right)$ and $g$-prime, then
\begin{equation}
R^{*}(a,b,p;q)=\frac{f(-q)}{f(-q^g)}
\end{equation}
\begin{equation}
\frac{d}{dq}R^{*}(a,b,p;q)=\frac{q^{-1}}{24}\frac{f(-q)}{f(-q^g)}(1-g-L(q)+gL(q^g)).
\end{equation}
\textbf{Proof.}\\
The proof follows from the above Lemma 2. 
\[
\]
Here we must mention that when $q=e^{-\pi\sqrt{r}}$, $L(q)$ can evaluated by
$$L_1(q)=\frac{1}{24}\left(1-(1-2x)z^2-6x(1-x)z\frac{dz}{dx}\right)=$$
$$=1-\frac{3\pi}{2\sqrt{r}}+K(k_r)\left(-1+2k_r^2+6\left(\frac{\alpha(r)}{\sqrt{r}}-k_r^2\right)K(k_r)\right) , $$
where
$$z=z(w)=K(w)={}_2F_{1}\left(\frac{1}{2},\frac{1}{2};1;w\right) , $$
$x=k_r$ is the modulus of $K(w)$ and $\alpha(r)$ is the elliptic alpha function (see [7] and [15]).
\[
\] 
\textbf{Lemma 4.} If $g$ is a prime number, then
$$
\left(\frac{n}{g^2}\right)\stackrel{?}{=}\sum^{\frac{g-1}{2}}_{j=1}X(j,0,g,n)+\frac{g-1}{2}X_g(n)
$$
and 
$$
\prod^{\frac{g-1}{2}}_{j=1}\left[j,g;q\right]_{\infty}=\frac{f(-q)}{f(-q^g)}
$$
\textbf{Proof.}
$$
R^{*}(a,b,p;q)=\prod^{\infty}_{n=1}(1-q^n)^{\left(\frac{n}{g^2}\right)}=f(-q^g)^{\frac{g-1}{2}}\prod^{\frac{g-1}{2}}_{j=1}R^{*}(j,0,g;q)=
$$
$$
=f(-q^g)^{\frac{g-1}{2}}\prod^{\frac{g-1}{2}}_{j=1}\frac{[j,g;q]_{\infty}}{f(-q^g)}=\prod^{\frac{g-1}{2}}_{j=1}[j,g;q]_{\infty}
$$
But also we can write 
$$R^{*}(a,b,p;q)=\prod^{\infty}_{n=1}(1-q^n)^{\left(\frac{n}{g^2}\right)}=\prod^{\infty}_{n=1}(1-q^n)\prod^{\infty}_{n=1}(1-q^{ng})^{-1}=f(-q)f(-q^g)^{-1}$$
Combining the above two results we get the proof.
\[
\]
\centerline{\textbf{Trying to Generalize the Problem}}
\[
\]
From the above it became obvious, that we can write for any $G$-integer with $G>1$:
\begin{equation}
\left(\frac{n}{G^2}\right)=\sum^{\nu}_{j=1}X(c_j,0,G;n)+cX_G(n) ,
\end{equation}
where $c$-integer constant depending from $G$, and $c_j$  positive integers depending also form $G$. 
\[
\]
\textbf{Theorem 12.} 
In the case which 
\begin{equation}
X(n)=\left(\frac{n}{G^2}\right)
\end{equation}
then
if $g_1,g_3,\ldots,g_{\lambda}$ are positive primes and $G=g_1^{w_1}g_2^{w_2}\ldots g_{\lambda}^{w_\lambda}$ is the factorization of $G$, with $g_1<g_2<\ldots<g_{\lambda}$, then  
$$
\prod^{\nu}_{j=1}[c_j,G;q]^{*}_{\infty}\stackrel{?}{=}f(-q)\prod^{\lambda}_{i=1}f(-q^{g_i})^{-1}\prod_{i<j}f(-q^{g_ig_j})^{1}\prod_{i<j<k}f(-q^{g_ig_jg_k})^{-1}\times $$
\begin{equation}
\times \prod_{i<j<k<l}f(-q^{g_ig_jg_kg_l})^{1}\ldots ,
\end{equation}
\begin{equation}
R^{*}(a,b,p;q)=\prod^{\nu}_{j=1}[c_j,G;q]^{*}_{\infty}
\end{equation}
The $c_j$ are given from (95).\\ Also
$$
\frac{d}{dq}R^{*}(a,b,p;q)=\frac{d}{dq}\prod^{\infty}_{n=1}(1-q^n)^{\left(\frac{n}{G^2}\right)}=
$$
$$
=q^{-1}f(-q)\left[L_1(q)-\sum^{\lambda}_{i=1}g_iL_1(q^{g_i})+\sum_{i<j}g_ig_jL_1(q^{g_ig_j})-\ldots\right]\times
$$
\begin{equation}
\times\prod^{\lambda}_{i=1} f(-q^{g_i})^{-1}\prod^{}_{i<j}f(-q^{g_ig_j})^{1}\prod^{}_{i<j<k}f(-q^{g_ig_jg_k})^{-1}\prod^{}_{i<j<k<l}f(-q^{g_ig_jg_kg_l})^{1}\ldots , 
\end{equation}
\[
\]
\textbf{Example.}\\
From the relation $$\left(\frac{n}{14}\right)^2=X(1,0,14;n)+X(3,0,14;n)+X(5,0,14;n)+3X_{14}(n)$$
we get
$$
\left[1,14;q\right]^{*}_{\infty}\left[3,14;q\right]^{*}_{\infty}\left[5,14;q\right]^{*}_{\infty}=\frac{f(-q)f(-q^{14})}{f(-q^2)f(-q^7)}
$$
\[
\]
\textbf{Note.}\\
If $X(a,b,p,n)$ can be written in the form $\left(\frac{n}{G}\right)$ with integer $G>1$, then
$$X(a,b,p;n)=\left(\frac{n}{G}\right)\stackrel{?}{=}\sum^{\left[\frac{G-1}{2}\right]}_{j=1}\left(\frac{j}{G}\right)X(j,0,G;n)$$ and 
$$
R^{*}(a,b,p;q)=\prod^{\infty}_{n=1}(1-q^n)^{\left(\frac{n}{G}\right)}=\prod^{\left[\frac{G-1}{2}\right]}_{j=1}\left([j,G;q]^{*}_{\infty}\right)^{\left(\frac{j}{G}\right)}
$$
Also 
$$
\prod^{G-1}_{j=1}[j,2G;q]^{*}_{\infty}=\frac{f(-q)}{f(-q^G)}
$$
\[
\]
\centerline{\textbf{The Characteristic Values Method}}
\[
\]
\textbf{Lemma 5.}\\ 
\textbf{i)} For every positive integers $a$, $b$, $G$ with $a<b<G$, exists integers $c(j)$, such that
\begin{equation}
X(a,b,G,n)=\sum^{G}_{j=1}c(j)X(j,0,G,n)
\end{equation}
\textbf{ii)} Suppose that $a$ and $G$ are integers. If exist integers $b(j)$, such that
\begin{equation}
X(a,0,G,n)=\sum^{G}_{j=1}b(j)X_j(n) , 
\end{equation}
then $R(a,b,G,q)$ can be determined by the functions $f(-q^g)^c$
\[
\]
\textbf{Example.}\\
For $a=2$, $b=6$, $p=G=12$, then
$$
X(2,6,12;n)=-X(6,0,12;n)+X(10,0,12;n)
$$ 
and
$$
X(6,0,12;n)=X_6(x)-2X_{12}(n)
$$
$$
X(10,0,12,n)=X_2(n)-X_4(n)-X_6(n)
$$
Hence
$$
X(2,6,12;n)=X_2(n)-X_4(n)-2X_{6}(n)+2 X_{12}(n)
$$
and
$$
R(2,6,12;q)=q^{2/3}\frac{f(-q^2)f(-q^{12})^2}{f(-q^4)f(-q^6)^2}$$

\section{Results and applications}

If  $v=R(1,2,5;q)=R(1,3,5;q)=R(q)$ is the Rogers Ramanujan continued fraction and $u=u(q)=R(1,3,10;q)$, then (see Routine 2, page 9):  
\begin{equation}
u^3-uv+u^2v^3+v^4\stackrel{?}{=}0
\end{equation}
In general its an observation that the Ramanujan Quantities are related with polynomial relations.\\
Another example is (see [10]):
\begin{equation}
k_r^2=16H^2\left(\frac{1-H^2}{1+H^2}\right)^2=\frac{(1-T)(3+T)^3}{(1+T)(3-T)^3} , 
\end{equation}
where $T=\sqrt{1-8V^3}$, $V=V(q)=R(1,3,6;q)$ is the Cubic continued fraction and $H=H(q)=R(1,3,8;q)$ is the Ramanujan-Gollnitz-Gordon continued fraction, $k_r$ is the elliptic singular moduli.\\
The $u$-quantity satisfies Theorem 8 hence 
\begin{equation}
u\left|u'\right|\stackrel{?}{=}u_2=u(q^2)  
\end{equation}
where $u'=u(-q)$\\
Another evaluation, using Mathematica is
$$
-v+(-1)^{1/5}v'-(-1)^{1/5}v^5v'+5(-1)^{2/5}v^4v'^2-10 (-1)^{3/5}v^3v'^3+
$$
\begin{equation}
+5(-1)^{4/5}v^2v'^4+vv'^5+v^6v'^5-(-1)^{1/5}v^5v'^6\stackrel{?}{=}0
\end{equation}
Relation (105) follows with the same method as (102) (see Routine 2 page 9)  
and $v'=|R(-q)|=v(-q)$.
\[
\]
\textbf{Proposition 3.} 
If $R(q)$ is the Rogers Ramanujan continued fraction, then 
\begin{equation}
R(q)R(q^2)\stackrel{?}{=}R(1,3,10;q)=u
\end{equation}
\textbf{Proof.}\\ 
It holds
\begin{equation}
\log(R(a,b,p;q))=Q\log(q)-\sum^{\infty}_{n=1}\frac{q^n}{n}\sum_{d|n}X(d)d , 
\end{equation} 
Multiplying with $\epsilon_n$ and summing we get  
\begin{equation}
\sum^{\infty}_{n=1}\frac{\epsilon(n)}{n}\log(R(a,b,p;q))=Q\log(q)\sum^{\infty}_{n=1}\epsilon(n)-\sum^{\infty}_{n=1}\frac{q^n}{n}\sum_{d|n}\tau(d)\epsilon(n/d)
\end{equation}
We make use of the properties of the Dirichlet Multiplication (see [3]).\\ Assume that exists arithmetic function $X_1(n)n$ such that 
$$X_1(n)n\ast 1=\tau(n)\ast \epsilon(n)=X(n)n\ast 1\ast \epsilon(n) , $$
then 
$$X_1(n)n=\epsilon(n)\ast X(n)n$$ 
Assuming the $X(n)$ and $X_1(n)$ are correspond to Ramanujan Quantities then we seek $\epsilon(n)$ such that is always 0 after a few terms. An example is when $X(n)=X(1,3,5,n)$ and $X_1(n)=X(1,3,10,n)$. Then 
$$
\epsilon(n)\stackrel{?}{=}1\textrm{ , if }n=1,2\textrm{ and } \epsilon(n)\stackrel{?}{=}0\textrm{ , }n>2
$$
Assuming the above and using (108) we have the proof of (106). 
\[
\]
\textbf{Note.}\\
1) With the method of Proposition 3 if we set  
$$
\log(R_{\lambda}(a,b,p;q)):=Q\log(q)-\sum^{\infty}_{n=1}\frac{q^n}{n}\sum_{d|n}X(a,b,p,\lambda d)d
$$
we have the evaluations 
$$R_3(1,3,8;q)=R_{11}(1,3,8;q)=qH(q)^{-1}$$
$$R_3(1,3,10;q)=q^{6/5}R(q)^{-1}R(q^2)^{-1}$$
$$
R_5(1,5,8;q)=qR(1,5,8;q)^{-1}
$$
$$
R_5(1,5,12;q)=q^{3/2}R(1,5,12;q)^{-1}
$$
$$R_{47}(1,5,7;q)=q^{5/7}R(1,3,7;q)^{-1}$$
\[
\]
\textbf{Application of Proposition 3.}\\
i) The function $m=R(q)R^2(q^2)$ can always be calculated in radicals with respect to $v=R(q)$, the equation that relates the two functions is
\begin{equation}
(m-1)\sqrt{vm}+v^3(1+m)=0
\end{equation}
or better
\begin{equation}
(m-1)\sqrt{R(q)m}+R(q)^3(1+m)=0
\end{equation}
ii) Moreover the fraction $v'=\left|R(-q)\right|$ can always evaluated in radicals with respect to $R(q)$.\\ 
\textbf{Proof.}\\ 
The proof of (i) follows from  
$u^2=R(q)R(q)R^2(q^2)$, hence $k=u^2/R(q)=u^2/v$, or $u=\sqrt{mv}$, setting this identity to (104) we get the result.\\
The proof of (ii) follows also from (104): $$u(q^2)=u(q)|u(-q)|=R(q)R(q^2)|R(-q)| R(q^2)=|R(-q)|m$$
Suppose the solution of (104) is $u(q)=f_1(R(q))$, then $u(q^2)=f_1(R(q^2))=|R(-q)|m$, hence 
$$|R(-q)|=\frac{f_1(R(q^2))}{m} . $$
Note that if $y=R(q^2)$, then
$$\frac{y-v^2}{y+v^2}=vy^2$$
The proof is complete.  
\[
\]
\centerline{\textbf{Applications of Theorem 12.}}
\textbf{1)}
$$
R(1,3,10;q)=R(q)R(q^2)=q^{3/5}\frac{f(-q^2)f(-q^{5})}{f(-q)f^3(-q^{10})}\left(\sum^{\infty}_{n=-\infty}(-1)^nq^{5n^2+4n}\right)^{2}
$$
\textbf{Proof.}\\
In [14] we have shown that 
$$
[a,p;q]^{*}_{\infty}=\frac{1}{f(-q^p)}\sum^{\infty}_{n=-\infty}(-1)^nq^{pn^2/2+(p-2a)n/2}
$$
Using the identity 
\begin{equation}
R(1,3,10;q)=q^{3/5}\frac{[1,10;q]^{*}_{\infty}}{[3,10;q]^{*}_{\infty}}=R(q)R(q^2)
\end{equation} 
and
\begin{equation}
[1,10;q]^{*}_{\infty}[3,10;q]^{*}_{\infty}=\frac{f(-q)f(-q^{10})}{f(-q^2)f(-q^5)}
\end{equation}
we get the result. 
\[
\]
\textbf{2)} Set $Y(q):=R(1,2,4;q)$,  $q=e^{-\pi\sqrt{r}}$, then
\begin{equation}
Y(q^4)=\sqrt{\frac{1-k'^2_r+2 \left(\sqrt{1+k'_r}-\sqrt{2}\right)t_1}{2 (1-k'_r)^2}}
\end{equation}
where
\begin{equation}
t_1=\sqrt{3 k'_r+k'^2_r+2\sqrt{2} k'_r\sqrt{1+k'_r}} .
\end{equation}
\textbf{Proof.}\\
If $R(1,2,4;q^8)=u=Y(q^8)$ and $R(1,3,8;q^2)=H(q^2)=v$ we have $$
-u+v+2u^2v-uv^2\stackrel{?}{=}0 .
$$
In view of article [10] and Theorem 10 of the present paper and relation 103, we get after simplifications the result.
\[
\]
\textbf{3)} Set $V_1(q)=R(1,3,12;q)$, then
\begin{equation}
V_1(q)=\frac{-1+T^{3/2}-\sqrt{3-2T^{3/2}-T^3}}{2(1-T^3)^{1/3}}
\end{equation}
\textbf{Proof.}\\
If $R(1,3,6;q^3)=V(q^3)=u$, where $V(q)$ is the Ramanujan's Cubic continued fraction (see [10]) and if $R(1,3,12;q^3)=v$ then
$$
u^3-uv+v^3+uv^4\stackrel{?}{=}0
$$
Hence in view of (103) with algebraic simplifications, we can always evaluate $V_1(q)=R(1,3,12;q)$, in radicals form, in terms of $T$ and hence of $k_r$ (for $k_r$ see also (103) and related references).
\[
\]
Other interesting results are:
\[
\] 
\textbf{a.} If $u:=R(1,3,6;q^3)=V(q^3)$,  
($V(q)$ is the Ramanujan's Cubic continued fraction, see [10]) and if $v:=R(1,3,12;q^3)$ then
$$
u^8+u^{11}v-u^7v^2-u^3v^3-u^{10}v^3+u^6v^4+3 u^2v^5-u^9 v^5+5u^5v^6-3uv^7+$$
$$+u^8 v^7+2 u^4 v^8+v^9\stackrel{?}{=}0
$$
\[
\]
\textbf{b.} 
If $u:=R(1,3,7;q)$, $v:=R(2,3,7;q)$, $w:=R(1,2,7;q)$ then
$$
u+u^3v-v^3\stackrel{?}{=}0\textrm{ and }-u^2+u^3w^2+w^3\stackrel{?}{=}0
$$
\[
\]
\textbf{4)}
\begin{equation}
R(1,2,6;q)=q^{1/4}\frac{f(-q)f^2(-q^6)}{f^2(-q^2)f(-q^3)}
\end{equation}
\textbf{Proof.}\\
It is easy to prove someone that
$$
[1,3;q]^{*}_{\infty}=\frac{f(-q)}{f(-q^3)}\textrm{ and }[1,6;q]^{*}_{\infty}=\frac{f(-q)f(-q^6)}{f(-q^2)f(-q^3)}
$$
Hence we get the result.
\[
\]
\textbf{5)} 
$$
\frac{(1+q^4)}{1+}\frac{q^8}{1+}\frac{(q^4+q^{12})}{1+}\frac{q^{16}}{1+}\frac{(q^{12}+q^{20})}{1+}\ldots
=\sqrt{2}\frac{q^{1/2}(1-k'_r)}{\sqrt{1-k'^2_r+2 \left(\sqrt{1+k'_r}-\sqrt{2}\right)t_1}}
$$
where
$$
t_1=\sqrt{3 k'_r+k'^2_r+2\sqrt{2} k'_r\sqrt{1+k'_r}} .
$$
\textbf{Proof.}\\ 
If $q=e^{-\pi\sqrt{r}}$, $r$ positive real
$\theta_4(q)=\sum^{\infty}_{n=-\infty}(-1)^nq^{n^2}$ from 
$$[1,2;q]^{*}_{\infty}=\prod^{\infty}_{n=0}(1-q^{2n+1})^2=X(-q)^2=\left(\frac{\theta_4(q)}{f(-q^2)}\right)^2$$
we get
$$[2,4;q]^{*}_{\infty}=X^2(-q^2)$$
and from $$[1,4;q]^{*}_{\infty}=\frac{f(-q)}{f(-q^2)}$$
we get (see and [17] introduction):
\begin{equation}
\frac{q^{1/8}}{Y(q)}=\frac{(1+q)}{1+}\frac{q^2}{1+}\frac{(q+q^3)}{1+}\frac{q^4}{1+}\frac{(q^3+q^5)}{1+}\ldots
=\frac{f(-q^2)X(-q^2)^2}{f(-q)} .
\end{equation}  
Using the above 2nd Application of Theorem 12 we get the result.
\[
\]
From
$$
[1,8;q]^{*}_{\infty}[3,8;q]^{*}_{\infty}=\frac{f(-q)}{f(-q^2)}
$$
we get
$$
R(1,2,8;q)=q^{5/16}\frac{f(-q)f(-q^4)}{f^2(-q^2)[3,8;q]^{*}_{\infty}}=
$$
$$=q^{5/16}\frac{f(-q)f(-q^4)f(-q^8)}{f^2(-q^2)}\left(\sum^{\infty}_{n=-\infty}(-1)^nq^{4n^2+n}\right)^{-1}
$$
Also hold
\begin{equation}
R(1,3,6;q)=V(q)=q^{1/3}\frac{[1,6;q]^{*}_{\infty}}{[3,6;q]^{*}_{\infty}}=q^{1/3}\frac{f(-q)f(-q^6)}{f(-q^2)f(-q^3)X^2(-q^3)}
\end{equation}
\textbf{6)}
\begin{equation}
R(2,3,6;q)=q^{1/12}\frac{f(-q^2)}{f(-q^6)X(-q^3)^2}
\end{equation}
\[
\]
\textbf{7)}
\begin{equation}
R^2(1,2,10;q)=\frac{f(-q)f^2(-q^{10})R(q)}{f^2(-q^2)f(-q^5)}
\end{equation}
\[
\]
\begin{equation}
R^2(1,4,10;q)=\frac{f(-q)f^2(-q^{10})R(q)R^2(q^2)}{f^2(-q^2)f(-q^{5})}
\end{equation}
\begin{equation}
R^2(2,3,10;q)=\frac{f^2(-q^2)f(-q^{5})R(q)R^2(q^2)}{f(-q)f^2(-q^{10})}
\end{equation}
\begin{equation}
R^2(3,4,10;q)=\frac{f(-q)f^2(-q^{10})}{f^2(-q^2)f(-q^{5})R(q)}
\end{equation}
\textbf{Proof.}\\
From the equations
$$[1,5;q]^{*}_{\infty}[3,5;q]^{*}_{\infty}=\frac{f(-q)}{f(-q^5)}$$
$$
[1,10;q]^{*}_{\infty}[3,10;q]^{*}_{\infty}=\frac{f(-q)f(-q^{10})}{f(-q^2)f(-q^5)}
$$
$$
R(1,3,5;q)R(1,3,5;q^2)=R(q)R(q^2)=R(1,3,10;q)
$$
and the definition of Ramanujan Quantities we get after solving the system the desired results.
\[
\]
Also from the above we have 
\begin{equation}
\left[1,10;q\right]^{*}_{\infty}=q^{-3/10} \sqrt{R(q)R(q^2)\frac{f(-q)f(-q^{10})}{f(-q^2)f(-q^5)}}
\end{equation}
and
\begin{equation}
\left[3,10;q\right]^{*}_{\infty}=q^{3/10}\sqrt{\left(R(q)R(q^2)\right)^{-1}\frac{f(-q)f(-q^{10})}{f(-q^2)f(-q^5)}}
\end{equation}
Hence one application is
\begin{equation}
R(1,3,10;q)=q^{4/5}\frac{\left[1,10;q\right]^{*}_{\infty}}{\left[3,10;q\right]^{*}_{\infty}}= \frac{q^{1/2}}{X(-q^5)^2}\sqrt{R(q)R(q^2)\frac{f(-q)f(-q^{10})}{f(-q^2)f(-q^5)}}  
\end{equation}
and because 1,5 are odd and 10 is even we have $$R(1,3,10;q^2)=|R(1,3,10;-q)|R(1,3,10;q) . $$ 
\[
\]
Last observe that $R(q)$ (Rogers-Ramanujan Continued fraction) can be written in terms of $f(-q^{\alpha})$ and $q^{\beta}$, using the following Ramanujan identities:
\begin{equation}
\frac{1}{R(q)}-1-R(q)=\frac{f(-q^{1/5})}{q^{1/5}f(-q^5)}
\end{equation}  
\begin{equation}
\frac{1}{R^5(q)}-11-R^5(q)=\frac{f^6(-q)}{qf^6(-q^5)} .
\end{equation}
Also observe that (128) is a modular equation of the agiles $x=[1,5;q]_{\infty}$ and $y=[3,5;q]_{\infty}$:
\begin{equation}
x^{10}-y^{10}+11x^5y^5+x^{11}y^{11}=0  
\end{equation}
Numerical calculations show us that such identities exists and for other 'agiles', examples are: 
\[
\]
i) If $x=[1,6;q^{3}]_{\infty}$ and $y=[3,6;q]_{\infty}$, then 
\begin{equation}
8x^9-y^3+x^{12}y^3+x^3y^6\stackrel{?}{=}0
\end{equation}
ii) If $x=[1,6;q^{2}]_{\infty}$ and $y=[2,6;q]_{\infty}$, then 
\begin{equation}
-9x^8+y^4+x^{12}y^4-x^4y^8\stackrel{?}{=}0
\end{equation}
iii) If $x=[1,4;q]_{\infty}$ and $y=[2,4;q]_{\infty}$, then 
\begin{equation}
16x^8+x^{16}y^4-y^8\stackrel{?}{=}0
\end{equation}
or
\begin{equation}
\frac{1}{R(1,2,4;q)^{16}}-\frac{16}{R(1,2,4;q)^8}=q^{-2}X(-q^2)^{24}
\end{equation}
iv) If $x=[1,6;q]_{\infty}$ and $y=[3,6;q]_{\infty}$, then 
\begin{equation}
8x^3-y^3+x^{12}y^3+x^9y^6\stackrel{?}{=}0
\end{equation}
Equation (134) lead us to the evaluation 
\begin{equation}
\frac{1-8V(q)^3}{V(q)^9(1+V(q)^3)}=q^{-3} X(-q^3)^{24}
\end{equation}
where $V(q)=R(1,3,6;q)$ is the Ramanujan's Cubic continued fraction.  
\[
\]
\newpage
\centerline{\bf References}\vskip .2in

[1]: I.J.Zucker, 'The summation of series of hyperbolic functions'. SIAM J. Math. Ana.10.192, 1979.

[2]: G.E.Andrews, Number Theory. Dover Publications, New York, 1994.

[3]: T.Apostol, Introduction to Analytic Number Theory. Springer Verlang, New York, Berlin, Heidelberg, Tokyo, 1974.

[4]: M.Abramowitz and I.A.Stegun, Handbook of Mathematical Functions. Dover Publications, 1972.

[5]: B.C.Berndt, Ramanujan`s Notebooks Part I. Springer Verlang, New York, 1985.

[6]: B.C.Berndt, Ramanujan`s Notebooks Part II. Springer Verlang, New York, 1989.

[7]: B.C.Berndt, Ramanujan's Notebooks Part III. Springer Verlang, New York, 1991.

[8]: L.Lorentzen and H.Waadeland, Continued Fractions with Applications. Elsevier Science Publishers B.V., North Holland, 1992. 

[9]: E.T.Whittaker and G.N.Watson, A course on Modern Analysis. Cambridge U.P., 1927.

[10]: Nikos Bagis, 'The complete evaluation of Rogers Ramanujan and other continued fractions with elliptic functions'. arXiv:1008.1304v1 [math.GM], 2010.

[11]: Michael Trott, 'Modular Equations of the Rogers-Ramanujan Continued Fraction'. page stored in the Web.

[12]: Nikos Bagis, 'The First Derivative of Ramanujans Cubic Fraction'. arXiv:1103.5346v1 [math.GM], 2011.  

[13]: Nikos Bagis. 'Parametric Evaluations of the Rogers-Ramanujan Continued Fraction'. IJMMS, Vol. 2011.

[14]: Nikos Bagis and M.L. Glasser. 'Jacobian Elliptic Functions, Continued fractions and Ramanujan Quantities'. arXiv:1001.2660v1 [math.GM], 2010.

[15]: J.M. Borwein and P.B. Borwein. 'Pi and the AGM'. John Wiley and Sons, Inc. New York, Chichester, Brisbane, Toronto, Singapore, 1987.

[16]: Nikos Bagis and M.L. Glasser. 'Integrals related with Rogers Ramanujan continued fraction and q-products'. arXiv:0904.1641v1 [math.NT], 2009.

[17]: Bhaskar Srivastava. 'Some Eisenstein Series Identities Related to Modular Equation of the Fourth Order'. Commun. Korean Math. Soc. 26 (2011), No. 1, pp. 67-77.

\end{document}